\documentclass{article}
\usepackage{amsmath}
\usepackage{amssymb}
\usepackage{latexsym}
\usepackage{amsthm}
\usepackage{mathrsfs} 
\usepackage{wasysym}
\usepackage{fancyhdr}
\usepackage{amsfonts}
\usepackage{amssymb}
\usepackage{epsfig}
\usepackage{epstopdf}
\usepackage[normalem]{ulem}
\usepackage{eucal}

\usepackage{tkz-graph}
\usetikzlibrary{arrows,decorations.markings}
\usepackage{tkz-euclide}
\usetkzobj{all}

\newtheorem{thm}{Theorem}[section]

\newtheorem{lemma}[thm]{Lemma}

\newtheorem{conj}[thm]{Conjecture}

\newtheorem{red}{Reduction}

\begin{document}
\title{Dominating maximal outerplane graphs and \\ Hamiltonian plane triangulations} 
\author{Michael D. Plummer\thanks{
Department of Mathematics,
Vanderbilt University,
Nashville, TN 37215,
Email: \tt{michael.d.plummer@vanderbilt.edu}}\,,
Dong Ye\thanks{Corresponding author. Department of Mathematical Sciences,
Middle Tennessee State University,
Murfreesboro, TN 37132, USA,
Email: {\tt{dong.ye@mtsu.edu}}. Partially supported by a grant from Simons Foundation (no. 359516).}\;  and
Xiaoya Zha\thanks{Department of Mathematical Sciences,
Middle Tennessee State University,
Murfreesboro, TN 37132, USA, Email: {\tt{xiaoya.zha@mtsu.edu}}.}}

\date{}
\maketitle
\begin{abstract}   
Let $G$ be a graph and $\gamma (G)$ denote the domination number of $G$, i.e. the cardinality of a smallest 
set of vertices $S$ such that every vertex of $G$ is either in $S$ or adjacent to a vertex in $S$.
Matheson and Tarjan conjectured that a plane triangulation with a sufficiently large number $n$ of vertices has $\gamma(G)\le n/4$.
Their conjecture remains unsettled. 
In the present paper, we show that: (1) a maximal outerplane graph with $n$ vertices has $\gamma(G)\le \lceil \frac{n+k} 4\rceil$
where $k$ is the number of pairs of consecutive degree 2 vertices separated by distance at least 3 on the boundary of $G$; and 
(2) a Hamiltonian plane triangulation $G$ with $n \ge 23$ vertices has $\gamma (G)\le  5n/16 $. 
We also point out and provide counterexamples for several recent published results of Li et al in [Discrete Appl. Math.198 (2016) 164-169] on this topic which are incorrect.

\medskip

\noindent {\it Keywords: plane triangulation, domination, Hamilton cycle, outerplane graph}
\end{abstract}

\section{  Introduction}%

A {\it plane triangulation} is a plane graph in which every face is bounded by a triangle. 
In 1996,  Matheson and Trajan made the following conjecture, in which the bound is tight if the conjecture is true. 

\begin{conj}[Matheson and Tarjan, \cite{MT}]
Let $G$ be plane triangulation with a sufficiently large number of vertices. Then $\gamma(G)\le |V(G)|/4$. 
\end{conj}

\begin{figure}[!hbtp] \refstepcounter{figure}\label{fig:2-graph}
\begin{center}
\includegraphics[scale=.75]{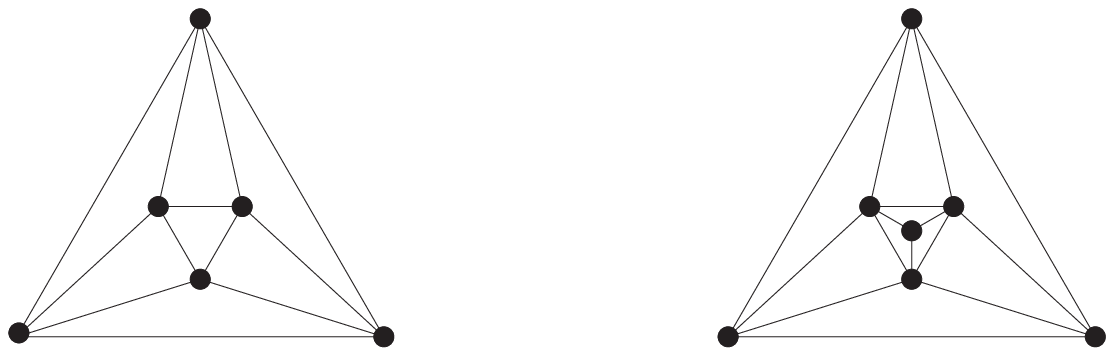} \\
{Figure \ref{fig:2-graph}: Two triangulations with large domination ratio.}
\end{center}
\end{figure}
 
A plane triangulation with small order $n$ may have domination number bigger than $n/4$. For example, the
triangle has $\gamma (K_3)=1 = n/3$, the octahedron shown in Figure~\ref{fig:2-graph} (left) has $\gamma=2=n/3$ and the 7-vertex graph shown in Figure~\ref{fig:2-graph} (right) has $\gamma =2= 2n/7>n/4$.  
So 
one must assume $n\ge 8$ in order for the Matheson-Tarjan conjecture to be true.
King and Pelsmajer \cite{KP} proved that the Matheson-Tarjan conjecture holds
when the maximum degree of the triangulation is 6.

Matheson and Tarjan \cite{MT} proved that 
every plane triangulation $G$ has domination number $\gamma(G)\le |V(G)|/3$,
and this bound was later proved for triangulations of the projective plane, torus and Klein bottle by 
Plummer and Zha \cite{PZ} 
and Honjo et al. \cite{HKN}.
Furuya and Matsumoto \cite{FM} generalized this result by showing that $\gamma (G)\le |V(G)|/3$,
for every triangulation $G$ of any closed surface.
More generally, the first and third authors of the present paper conjectured that if a
triangulation of {\it any}
non-spherical surface with $n$ vertices, then
$\gamma (G)\le n/4$ if $n$ is sufficiently large \cite{PZ}.
Both these conjectures involving the $n/4$ bound remain unsettled.

In \cite{PYZ}, the authors
of the present paper proved that a Hamiltonian plane triangulation with minimum degree at least 4 has domination number at most
$\max\{\lceil 2n/7\rceil, \lfloor 5n/16 \rfloor\}$.
Very recently,  \v{S}pacapan \cite{SP} show that every plane triangulation has $\gamma (G)\le 17n/53$. 
It turns out that the domination number of maximal outerplane graphs plays a very important role
in the proofs of both of the main results in \cite{PYZ}
and \cite{SP}.

An {\it outerplane graph} is a graph embedded in the plane in such a way that all vertices of the graph lie on the boundary of the infinite face.
An outerplane graph is {\it maximal} if it is not possible to add any new edge to $G$ without destroying outerplanarity.
So a maximal outerplane graph is a 2-connected outerplane graph in which every inner face is a triangle. 
In 2013, Campos and Wakabayashi \cite{CW} obtained the following result for maximal outerplane graphs. 

\begin{thm}[Campos and Wakabayashi \cite{CW}, Tokunaga\cite{ST}]\label{thm:outplanar}
If $G$ is a maximal outerplane graph with $n\ge 4$ vertices and with $t$ vertices of degree 2, then $\gamma (G)\le (n+t)/4$.
\end{thm}

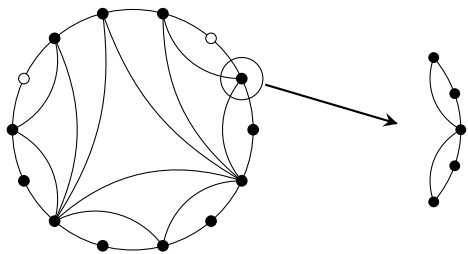
\begin{figure}[!hbtp]
\begin{center}
\begin{tikzpicture}[scale=.8]
 \tikzset{
  big arrow/.style={
    decoration={markings,mark=at position 1 with {\arrow[scale=1.3,#1]{>}}},
    postaction={decorate},
    shorten >=0.4pt},
  big arrow/.default=black}

\draw[] (0,0) circle (2);
\draw[] (.5,1.93) to [bend right=40]   (1.81,.85)   to [bend right=40]   (1.81,-.85)   to [bend right=40]   (.5,-1.93)    to [bend right=40]    (-1.3,-1.52)   to [bend right=40]  (-2,0)    to [bend right=40]   (-1.3, 1.52);

\draw[] (-.5,1.93)  to [bend right=20]  (1.81,-.85)   to [bend right=30] (-1.3,-1.52)  to [bend right=20] (-.5,1.93);

\draw[]  (-1.3,-1.52)   to [bend right=25]  (-1.3,1.52) ; \draw[] (.5,1.93) to [bend right=25]  (1.81,-.85);

\filldraw[] (-2,0)  circle (2.5pt);  \filldraw[] (2,0) circle(2.5pt); \filldraw[] (-.5,1.93) circle(2.5pt); \filldraw[] (.5,1.93) circle(2.5pt); 
\filldraw[] (-1.3, 1.52) circle(2.5pt);


\filldraw[fill=white, draw=black] (-1.81, .85) circle(2.5pt); 

\filldraw[fill=white, draw=black] (1.3, 1.52) circle(2.5pt); 

\filldraw[] (1.81,.85) circle(2.5pt);

\filldraw[] (-.5,-1.93) circle(2.5pt); \filldraw[] (.5,-1.93) circle(2.5pt); 
\filldraw[] (-1.3,-1.52) circle(2.5pt); \filldraw[] (-1.81,-.85) circle(2.5pt);\filldraw[] (1.3,-1.52) circle(2.5pt); \filldraw[] (1.81,-.85) circle(2.5pt);
 
 \draw[thick, big arrow,>=stealth] (2.2,.75)--(4.4,0.1); 
 
 \draw[] (1.81,.85) circle (10pt);
 
 \draw[] (5,1.2)  to [bend left=40]  (5,-1.2);
 \draw[] (5, 1.2) to [bend right=40] (5.45, 0) to [bend right=40] (5,-1.2);
  \filldraw[] (5,1.2) circle(2.3pt); \filldraw[] (5.45,0) circle(2.3pt);  \filldraw[] (5,-1.2) circle(2.3pt); 
  \filldraw[] (5.35,.6) circle(2.3pt);  \filldraw[] (5.35,-.6) circle(2.3pt); 

 \end{tikzpicture}
\caption{\small A maximal outerplane graph with domination number 4.}\label{fig:counterexample}
\end{center}
\end{figure}

Let $G$ be a maximal outerplane graph. Then the outer boundary of $G$ is a Hamilton cycle $C$ of $G$.
Let $C=v_1v_2\cdots v_nv_1$ be this Hamilton cycle such that $v_1, v_2, \ldots, v_n$ appear in clockwise
order along the boundary of $G$.
Two degree 2 vertices $v_r$ and $v_s$ are
{\it consecutive}
if the segment $C[r+1,s-1]=v_{r+1}\cdots v_{s-1}$ does not contain any degree 2 vertices.
Let an ordered pair of consecutive degree 2 vertices $(v_r, v_s)$ be called
{\em essential} if $|E(C[r,s])|\ge 3$. (For example, the pair of white vertices in Figure~\ref{fig:counterexample} is an essential pair.)
The first vertex in an essential pair is called a {\em bad} vertex in \cite{LZSX}.
Clearly, the number of essential pairs of $G$ is equal to the number of bad vertices of $G$.   

In 2016, Li, et al. \cite{LZSX}  realized that the degree $2$ vertices which require more vertices to dominate are
precisely the bad vertices.
They claimed that the domination number for an $n$-vertex maximal outerplane graph with $k$ bad vertices (or essential pairs)
is at most $(n+k) / 4$. 
However, this upper bound is not correct. 
For example, the graph $G$ in Figure~\ref{fig:counterexample} has $n=14$ and $k=1$, but $\gamma(G)=4> (n+k)/4$.
It is easy to extend the graph $G$ in Figure~\ref{fig:counterexample} to an infinite family of maximal outerplane graphs
which do not satisfy the bound of Li et al. by replacing a non-degree 2 vertex of $G$ by a $4m+1$ path with $2m$ consecutive triangles
and then adding interior chords to form a maximal outerplane graph. (In Figure 2 we indicate this operation in
the case when m = 1.) The smallest counterexample to the incorrect bound is a 6-vertex maximal outerplane graph $G$ with exactly
three degree 2 vertices. 

In this paper, we first obtain a bound for the domination number of a maximal outerplane graph based on essential pairs,
which corrects the result of Li et al. \cite{LZSX}, as follows.


\begin{thm}\label{thm:main1}
If $G$ is a maximal outerplane graph with $n\ge 4$ vertices and $k$ essential pairs, then $\gamma (G)\le \lceil (n+k)/4 \rceil$.
\end{thm}

Another main result (Theorem 3.2) of Li et al. \cite{LZSX} claims that every Hamiltonian plane triangulation
with at least seven vertices has domination number at least $5n/16$, which is better than the result obtained in \cite{PYZ}.
However, the proof of Theorem 3.2 of \cite{LZSX} heavily depends on the incorrect bound $\gamma(G)\le (n+k)/4$
for maximal outerplane graphs, and Lemma~3.1 of \cite{LZSX}.
Unfortunately, the proof of Lemma~3.1 in \cite{LZSX}
contains errors.
For example, in the proof of their Claim, if $v_{i_l+1}$ is a degree 2 vertex of $G_{out}^C$ and $v_{i_l-2}$ is a degree 2 vertex
of $G_{in}^C$, the rerouting reduces the distance between $v_{i_l}$ and $v_{i_{l+1}}$,
but increases the distance between $v_{i_l-2}$ and $v_{i_l}$. If $v_{i_l-2}$ is a degree 2 vertex, then it is not bad in $G_{in}^C$, but is bad in $G_{in}^{C'}$. Therefore, the total number of bad vertices of $G_{in}^{C'}$ remains
unchanged after the rerouting.  
We do not see a way to
correct
these mistakes. It remains unknown whether the result of Lemma 3.1 of \cite{LZSX} is true or not. 

In the present paper, using a recent result of Brinkmann, Ozeki and Van Cleemput \cite{BOV},
we also are able to prove the following result, which strengthens the main result of \cite{PYZ}
and
is
almost as good as the result claimed by Li et al. in \cite{LZSX}.

\begin{thm} \label{thm:main2}
Every Hamiltonian plane triangulation with at least $n\ge 23$ vertices has domination number at most $5n/16$. 
\end{thm}

\section{Maximal outerplane graphs} 
In this section, we will prove
Theorem~\ref{thm:main1}.
Before proceeding to prove Theorem~\ref{thm:main1}, we first introduce some reductions.
In these reductions, a maximal outerplane graph $G$ is reduced to a smaller maximal outerplane graph $G'$.
For convenience, we always assume $|V(G)|=n$ and $|V(G')|=n'$, and denote the numbers of essential pairs of $G$ and $G'$ by $k$ and $k'$, respectively. 

Let $G$ be a maximal outerplane graph with $n$ vertices.
The boundary of $G$ is a Hamilton cycle of $G$, which
we will denote
by $C=v_1v_2\cdots v_nv_1$ in clockwise order around the boundary.
For any two vertices $v_r$ and $v_s$, denote by $C[r,s]$ the segment of $C$ from $v_r$ to $v_s$ in the clockwise direction. If $v_rv_s$ is a chord, the subgraph $G[r,s]$  induced by the vertices on
$C[r,s]$ is called a {\em section} of $G$.
All vertices of $C[r,s]$
different from
$v_r$ and $v_s$
will be
called {\it internal vertices} \rm of $C[r,s]$ (or $G[r,s]$).
Clearly, $G[r,s]$ contains degree $2$ vertices of $G$ as its internal vertices since $v_r v_s$ is a chord of $C$.
A section $G[r,s]$ is {\it elementary} if it contains exactly one degree 2 vertex of $G$.
An elementary section of $G$ is {\it maximal } if it is not contained in any other elementary section of $G$.
A triangular face $v_rv_sv_t$ of $G$ is called an {\em internal triangle} if all three edges of $v_rv_sv_t$ are chords of $C$.
An internal triangle $v_rv_sv_t$ of $G$ separates $G$ into three sections $G[r,s]$, $G[s,t]$ and $G[t,r]$.
A maximal outerplane graph is {\em striped} if it has no internal triangle.
In other words, the inner dual of a striped maximal outerplane graph is a path.
In the following reductions, the vertex label remains unchanged if the reduction does not affect the vertices.
Note that an essential pair $(v_r,v_s)$ 
of $G$
remains
an essential pair of $G'$ if the reduction does not affect the segment $C[r,s]$.

\medskip

\begin{red}\label{red1}   
Assume that $G[r,s]$ is an elementary section of $G$ with six vertices.  
Delete
the four
internal vertices of $C[r,s]$ to obtain a new maximal outerplane graph $G'$. 
\end{red}
 
\noindent{\bf Claim:}  {\sl $n'= n-4$, $k'\le k$ and $\gamma(G)\le \gamma(G')+1$. } \medskip

\noindent{\em Proof of Claim.} It is clear that $n'=n-4$. Let $v_{t_1}, v_{t_2}$ and $v_{t_3}$ be the three consecutive degree 2 vertices
on the boundary $C$ of $G$ such that $v_{t_2}\in C[r,s]\subset C[t_1,t_3]$.
Since $G$ is a maximal outerplane triangulation, at most one of $v_r$ and $v_s$ is a degree 2 vertex of $G'$.
First, assume that 
neither of them is a degree 2 vertex.
Note that at least one of $(v_{t_1}, v_{t_2})$ and $(v_{t_2},v_{t_3})$ is an essential pair of $G$ because $C[r,s]$ has
five
edges.
But $G'$ does not contain $v_{t_2}$ and $(v_{t_1}, v_{t_3})$ may be an essential pair of $G'$. Hence, $k'\le k$.
Now assume one of $v_r$ and $v_s$ is a degree 2 vertex.
By symmetry, assume $v_r$ is a degree 2 vertex of $G'$.
Let $C'$ be the boundary of $G'$.
Then $C'=(C-C[r,s])\cup \{v_rv_s\}$.
It follows that $|E(C[t_1, t_2])|\ge |E(C'[t_1,r])|$ and $|E(C[t_2, t_3])|\ge |E(C'[r,t_3])|$.
Therefore, $(v_{t_1},v_{t_2})$ is an essential pair of $G$ if $(v_{t_1}, v_r)$ is an essential pair of $G'$, and $(v_{t_2},v_{t_3})$ is an essential pair of $G$ if $(v_{r}, v_{t_3})$ is an essential pair of $G'$. So $k'\le k$ follows. 

Since $G[r,s]$ has six vertices, it follows that $s=r+5$ and hence $v_s=v_{r+5}$.
By symmetry, assume $v_{r+1}$ or $v_{r+2}$ is the degree $2$ vertex.
Let $D'$ be a dominating set of $G'$. Note that $G'$ is a subgraph of $G$. The vertices $v_r$ and $v_s$ are dominated by $D'$. 
If $v_{r+1} $ is the degree $2$ vertex, then $v_r$ is adjacent to
each of the other five vertices of $G[r,s]$,
and so $D' \cup v_r$ is a dominating set of $G$.
Now suppose $v_{r+2}$ is the degree $2$ vertex.
Then $D' \cup v_{r+3}$ is a dominating set of $G$.  
Hence $\gamma(G)\le \gamma(G')+1$ and the Claim is proved.
 \qed\medskip

 \begin{red}\label{red2}
Let $G[r,s]$ be a section with at least six vertices and suppose it contains a vertex adjacent to all other vertices of $G[r,s]$.
Delete all internal vertices of $G[r,s]$ to generate a new maximal outerplane graph $G'$.
 \end{red}

\noindent{\bf Claim:}  {\sl $n'\le n-4$, $k'\le k$ and $\gamma(G)\le \gamma(G')+1$. } \medskip

\noindent{\em Proof of Claim.} Since $G[r,s]$ has at least six vertices, it follows that  $n'\le n-4$.

Let $v_i$ be the vertex dominating all vertices of $G[r,s]$.  Given a dominating set $D'$ of $G'$, the set $D=D'\cup\{v_i\}$ is a dominating set of $G$.  Hence, $\gamma(G)\le \gamma(G')+1$. In the following, we show $k'\le k$.

If $v_i=v_r$ or $v_s$, then $G[r,s]$ has exactly one degree 2 vertex of $G$ which is either $v_{r+1}$ or $v_{s-1}$. An argument similar to that used in the proof of Reduction~\ref{red1} shows that $k\le k'$. 
So assume that $v_i$ is an internal vertex of $C[r,s]$. Then $G[r,s]$ has exactly two degree 2 vertices which are $v_{i-1}$ and $v_{i+1}$. 
Let $v_{t_1}$ and $v_{t_2}$ be two degree 2 vertices such that neither $G[t_1, r]$ nor $G[s,t_2]$ contains a degree 2 vertex
as an internal vertex. Since $G[r,s]$ has at least six vertices, it follows that either one of $(v_{t_1},v_{i-1})$ and $(v_{i+1}, v_{t_2})$ or both  are essential pairs. After the reduction, at most one of $v_r$ and $v_s$ is a new degree 2 vertex in $G'$.  Let $C'$ be the boundary of $G'$ which is $(C-C[r,s])\cup v_rv_s$. Then $|E(C[t_1, i-1])|\ge 1+ |E(C[t_1,r])|=|E(C'[t_1,s])|>|E(C'[t_1,r])|$, and $|E(C[i+1, t_2])|\ge 1+|E(C[s,t_2])|=|E(C'[r,t_2])|>|E(C'[s,t_2])|$. 

So, if $v_r$ is the new degree 2 vertex, then $(v_{t_1}, v_{i-1})$ is an essential pair of $G$ if $(v_{t_1},v_r)$ is an essential pair of $G'$, and $(v_{i+1}, v_{t_2})$ is an essential pair of $G$ if $(v_{r}, v_{t_2})$ an essential pair of $G'$; and
if $v_s$ is the new degree 2 vertex, then $(v_{t_1}, v_{i-1})$ is an essential pair of $G$ if $(v_{t_1},v_s)$ is an essential pair of $G'$, and $(v_{i+1}, v_{t_2})$ is an essential pair of $G$ if $(v_{s}, v_{t_2})$ is an essential pair of $G'$. Therefore, $k\le k'$ if  one of $v_r$ and $v_s$ is a degree 2 vertex of $G'$. 

So assume that neither $v_r$ nor $v_s$ is a degree 2 vertex, then $(v_{t_1},v_{t_2})$ is the only essential pair of $G'$, which is not an essential pair of $G$.  Hence $k'\le k$. This completes the proof that $k'\le k$.  \qed\medskip

\begin{red}\label{red3} 
Assume that $G[r,s]$ is an elementary section with five vertices such that the middle internal vertex is not a degree 2 vertex of $G$.  
Contract $G[r, s]$ to a single vertex $x$ and delete all resulting multiple edges to obtain a new maximal outerplane graph $G'$.
\end{red}

\noindent{\bf Claim:}  {\sl $n'=n-4$, $k'\le k$ and $\gamma(G)\le \gamma(G')+1$. } \medskip

\noindent{\em Proof of Claim.}  Since $G[r,s]$ has exactly five vertices, it follows that $s=r+4$ (i.e., $v_s=v_{r+4}$), and  $n'=n-4$.

Let $t_1$ and $t_2$ be two degree 2 vertices of $G$ such that both $G[t_1, r]$ and $G[s,t_2]$ do not contain a degree 2 vertex
as an internal vertex. Let $v_i$ be the degree 2 vertex of $G[r,s]$, which is an internal vertex of $C[r,s]$.
It follows that $|E(C'[t_1, x])|<|E(C[t_1, i])|$ and  $|E(C'[x, t_2])|< |E(C[i, t_2])|$. Hence $k'\le k$ if $x$ is a degree 2 vertex of $G'$. Hence, assume that  the  vertex $x$ is not a degree 2 vertex of $G'$. 
If $(t_1,t_2)$ is an essential pair of $G'$, then $|E(C[t_1, r])|+|E(C[s,t_2])|\ge 3$. So $|E(C[t_1,i])|+|E(C[i,t_2])|=|E(C[r,s])|+|E(C[t_1, r])|+|E(C[s,t_2])|\ge 8$. Therefore, at least one of $(t_1,v_i)$ and $(v_i, t_2)$ is an essential pair of $G$. Hence  $k'\le k$.

Let $D'$ be a dominating set of $G'$.
By the symmetry of $v_{r+1}$ and $v_{s-1}$, we may assume $v_{r+1}$  is the degree $2$ vertex.
Then $v_r$ is adjacent to all vertices of $G[r,s]$. 
If $x\in D'$, then $D= (D'-x)\cup \{v_r  ,v_s \}$ is a dominating set of $G$.
If $x \not \in D'$, then $D=D' \cup v_r$ is a dominating set of $G$. So $\gamma(G)\le \gamma(G')+1$. 
 \qed\medskip

\begin{red}\label{red4}
Assume $v_rv_sv_t$ spans an internal triangle such that both $G[r,s]$ and $G[s,t]$ have exactly one internal vertex.
Contract $G[r,t]$ to a single vertex $x$ and delete any resulting multiple edges to obtain a new maximal outerplane graph $G'$. 
\end{red} 

\noindent{\bf Claim:}  {\sl $n'=n-4$, $k'\le k$ and $\gamma(G)\le \gamma(G')+1$. } \medskip

\noindent{\em Proof of Claim.} The section $G[r,t]$ has exactly five vertices.
So $n'=n-4$.

Since each of $G[r,s]$ and $G[s,t]$ has one internal vertex, $v_{r+1}$ and $v_{t-1}$ are degree 2 vertices of $G$.
Let $v_{t_1}$ and $v_{t_2}$ be the two degree 2 vertices on the boundary $C$ of $G$ such neither $C[t_1,r]$ nor $C[t,t_2]$
contains a degree 2 vertex as an internal vertex.
Neither of $v_{r+1}$ and $v_{s+1}$ is a vertex of $G'$. 
Note that $|E(C'[t_1, x])|<|E(C[t_1, r+1])|$ and  $|E(C'[x, t_2])|< |E(C[t-1, t_2])|$. 
It follows that $k'\le k$ if $x$ is a degree 2 vertex of $G'$. So assume that 
$x$ is not a degree 2 vertex of $G'$.
If $(v_{t_1},v_{t_2})$ is an essential pair of $G'$, then either one of $C[t_1, r]$ and $C[t,t_2]$ or both have at least two edges.
Then either one of $(v_{t_1},v_{r+1})$ and $(v_{t-1},v_{t_2})$ or both are essential pairs of $G$.
Hence, $k'\le k$. 

Let $D'$ be a dominating set of $G'$.
If $x\in D'$, then $D=(D'-x)\cup \{v_r, v_t\}$ is a dominating set of $G$.
If $x\notin D'$, then $D=D'\cup \{v_s\}$ is a dominating set of $G$.
No matter which of these two cases occurs, it follows that $\gamma(G)\le \gamma(G')+1$. 
 \qed\medskip

Reductions \ref{red1} and \ref{red2} can be applied to any maximal plane graph as long as they contain such sections. However, it requires 
the maximal plane graph to have at least seven vertices so that Reductions \ref{red3} and \ref{red4} can be applied to generate a new maximal outerplane 
graph (with at least three vertices).  

Each of the Reductions~\ref{red1}--\ref{red4} can be used to obtain upper bounds for the domination number by preserving the $1/4$ ratio
(required by the Matheson-Tarjan conjecture) or any other
ratio $\alpha$ between $1/4$ and $1/3$.  
A maximal plane triangulation $G$ is {\em $\alpha$-reducible} if it can be reduced to a smaller maximal plane triangulation $G'$
using 
some 
reductions such that $\gamma(G)\le \lceil \alpha(n+k)\rceil$ if $\gamma(G')\le \lceil \alpha(n'+k')\rceil$.
Otherwise, $G$ is {\em $\alpha$-irreducible}. 
In all Reductions~\ref{red1}--\ref{red4},  if $\gamma(G')\le \lceil \alpha(n'+k')\rceil$ for $\alpha\ge 1/4$, then
\[\gamma(G)\le \gamma(G')+1\le \lceil \alpha(n'+k') \rceil +1\le  \lceil \alpha(n-4+k) \rceil +1\le  \lceil \alpha(n+k) \rceil,\]
and hence $G$ is $\alpha$-reducible.  So we have the following result for $\alpha$-irreducible maximal outerplane graphs. 

\begin{lemma}\label{lem:reduce}
Let $G$ be a maximal outerplane graph with at least seven vertices. If $G$ is $\alpha$-irreducible for $\alpha \ge 1/4$, then:\\  
(1) Every maximal elementary section has at most three internal vertices;
furthermore, if a maximal elementary section has exactly three internal vertices, then the degree 2 vertex must be the middle internal vertex.  \\
(2) Two consecutive sections $G[r,s]$ and $G[s,t]$ have a total of at least three internal vertices if $v_rv_t$ is a chord.\\
(3) Any section with at least six vertices has no dominating vertex. 
 \end{lemma}

We now proceed to prove Theorem~\ref{thm:main1}. \medskip

\noindent{\bf Proof of Theorem~\ref{thm:main1}.} It is easy to check that the theorem holds for all maximal outplane graphs with 
at most six vertices.  In the following, assume that $G$ is a minimum counterexample. Then $G$ is $\alpha$-irreducible for $\alpha\ge 1/4$ and has at least seven vertices. 
Let $C$ be the boundary of $G$ oriented clockwise.  
\medskip

\noindent {\bf Claim~1.} {\sl Every elementary section $G[r,s]$ has at most two internal vertices.} 
\medskip

\noindent{\em Proof of Claim~1.} Suppose to the contrary that $G[r,s]$ has at least three internal vertices.
By (1) in Lemma~\ref{lem:reduce}, $G[r,s]$ has exactly three internal vertices and the degree 2 vertex is $v_{r+2}=v_{s-2}$.
Let $G'=G-\{v_{r+1}, v_{r+2}, v_{r+3}\}$.
Then $G'$ is a smaller maximal outerplane graph that $G$.
Since $G$ is a minimum counterexample, the domination number of $G'$ satisfies $\gamma(G')\le \lceil (n'+k')/4\rceil$. 

Let $v_{t_1}$ and $v_{t_2}$ be two degree 2 vertices such that $v_{t_1}, v_{r+2}$ and $v_{t_2}$ are three consecutive degree 2 vertices of $G$ appearing
in the clockwise direction of $C$.
Then $C[r,s]\subset C[t_1,t_2]$ and hence both $(v_{t_1}, v_{r+2})$ and $(v_{s_2}, v_{t_2})$ are essential pairs of $G$.
However, $G'$ contains all essential pairs of $G$ except $(v_{t_1}, v_{r+2})$ and $(v_{s_2}, v_{t_2})$,
together with the new  essential pair $(v_{t_1},v_{t_2})$.
Hence, $k'=k-1$.
It then follows that $n'+k'=n+k-4$. 

Now, let $D'$ be a dominating set of $G'$.
If $v_r\in D'$, let $D=D'\cup \{v_{r+3}\}$.
If $v_r\notin D'$, let $D=D'\cup \{v_{r+1}\}$.
Then $D$ is a dominating set of $G$.
Hence $\gamma(G)\le \gamma(G')+1\le  \lceil (n'+k')/4\rceil +1 \le  \lceil (n+k-4)/4\rceil+1\le  \lceil (n+k)/4\rceil$,
which contradicts the assumption that $G$ is a counterexample. \medskip

If $G$ has no internal triangle, then $G$ is a striped maximal outerplane graph which has exactly two degree 2 vertices 
 $v_i$ and $v_j$.
(Note: see page~3 
for the definition of a striped maximal outerplane graph which was first defined in \cite{CW}.)
These two degree 2 vertices either have a common neighbor or $\min\{|E(C[i,j])|, |E(C[j,i])|\}\ge 3$.
If the former holds, $G$ has a dominating vertex, namely the common neighbor of $v_i$ and $v_j$.
If the latter holds, then the total number of degree 2 vertices of $G$ is equal to $k=2$, and
Theorem~\ref{thm:main1} follows from Theorem~\ref{thm:outplanar}.

So in the following we shall assume that $G$ has at least one internal triangle.
Then $G$ has an internal triangle $T=v_rv_sv_t$ where $v_r, v_s$ and $v_t$ appear in
the clockwise direction of $C$ such that $G[r,s]$ and $G[s,t]$ are two consecutive elementary sections.
Let $v_{t_1}$ and $v_{t_2}$ be the two degree 2 vertices such that $C[t_1,t_2]$ is the shortest segment of $C$ containing $C[s,t]$.
 \medskip
 
\noindent{\bf Claim~2.} {\sl At least one of $G[r,s]$ and $G[s,t]$ contains exactly one internal vertex.}
\medskip

\noindent{\em Proof of Claim~2.} Suppose to the contrary that both $G[r,s]$ and $G[s,t]$
contain at least two internal vertices.
By Claim~1, both of them then have {\it exactly} two internal vertices.
By (3) of Lemma~\ref{lem:reduce}, at most one of $v_{r+2}=v_{s-1}$ and $v_{s+1}=v_{t-2}$ is a degree 2 vertex;
otherwise, $v_s$ is a dominating vertex of the section $G[r,t]$.
Then, at least one of $v_r$ and $v_t$ is adjacent to a degree 2 vertex.

First, assume that both $v_{r+1}$ and $v_{t-1}$ are degree 2 vertices.
Then contract $C[r+1,t-1]$ to a new vertex $x$ and delete any resulting multiple edges to obtain a smaller maximal outerplane graph $G'$
with $n'=|V(G')|=n-4$.
Note that the new vertex $x$ has degree 2 and is adjacent to both $v_r$ and $v_t$ in $G'$.
Hence, $G'$ has one less essential pair than $G$ because $(v_{r+1},v_{t-1})$ is an essential pair of $G$, but not $G'$.
Since $G$ is a minimum counterexample, we have $\gamma(G')\le  \lceil (n'+k')/4\rceil$.
Let $D'$ be a dominating set of $G'$.
Then at least one of $x$, $v_{r}$ and $v_t$ belongs to $D'$ because $x$ is a degree 2 vertex.
Then the set $D=(D'-x)\cup \{v_r,v_t\}$ dominates all the vertices of $G$.
Note that $|D|\le |D'|+1$, because $\{v_r,v_t\}\cap D'\ne \emptyset$ if $x\notin D'$.
It follows that $\gamma(G)\le \gamma(G')+1\le  \lceil (n'+k')/4\rceil+1= \lceil (n-4+k-1)/4\rceil+1 \le  \lceil (n+k)/4\rceil$,
a contradiction of the assumption that $G$ is a counterexample.
So at most one of $v_{r+1}$ and $v_{t-1}$ is a degree 2 vertex. 

By the symmetry of $v_{r+1}$ and $v_{t-1}$ in $G[r,t]$, we may assume that $v_{r+1}$ is a degree 2 vertex.
Then $v_{t-1}$ is not a degree 2 vertex.
Instead of $v_{t-1}$, the vertex $v_{s+1}$ is a degree 2 vertex.
Now contract $G[s,t]$ to a single vertex $x$ to form a smaller maximal outerplane graph $G'$ which has $n'=|V(G')|=n-3$ vertices.
(If $v_{r+1}$ is not a degree 2 vertex, then contract $G[r,s]$ to a new single vertex $x$.)
Since $G$ is a minimum counterexample, it follows that $\gamma(G')\le \lceil (n'+k')/4\rceil$.
Then $|E(C'[r+1, t_2])|=|E(C[s+1, t_2])|$ and $C'[t_1, r+1]=C[t_1, r+1]$.
Therefore, $G'$ has one less essential pair than $G$ because $(v_r, v_s)$ is not an essential pair of $G'$ and
$(v_s, v_{t_2})$ is replaced by $(v_r, v_{t_2})$ in $G'$.
Hence $k'\le k-1$, and further $n'+k'\le n+k-4$. 
Let $D'$ be a dominating set of $G'$ and let $D=D'\cup \{v_s\}$.
Then $D$ is a dominating set of $G$.
Hence $\gamma(G)\le \gamma(G')+1\le  \lceil (n'+k')/4\rceil+1= \lceil (n-4+k-1)/4\rceil+1 \le  \lceil (n+k)/4\rceil$,
again contradicting that $G$ is a counterexample.
This completes the proof of Claim~2. 
\medskip

By Claim~1, Claim~2 and (2) in Lemma~\ref{lem:reduce}, one of $G[r,s]$ and $G[s,t]$ has one internal vertex and the other
has two internal vertices.
By (3) in Lemma~\ref{lem:reduce}, then $v_{r+1}$ and $v_{t-1}$ are degree 2 vertices of $G[r,t]$ which form an essential pair.
Note that $G[r,t]$ is dominated by $\{v_r,v_t\}$. 

Contract $G[r+1, t-1]$ to a new vertex $x$ which is then a degree 2 vertex in the newly formed maximal outerplane triangulation $G'$.
Note that, $(v_{t_1}, x)$ is an essential pair of $G'$ if $(v_{t_1}, v_r)$ is an essential pair of $G$, and the same holds for $(x,v_{t_2})$.
But, $(v_{r+1},v_{t-1})$ is not an essential pair of $G'$.
Hence, $G'$ has one essential pair less than $G$, which means $k'=k-1$.
Since $G'$ has fewer vertices than $G$, it follows that $\gamma(G')\le \lceil (n'+k')/4\rceil$.
For every dominating set $D'$ of $G'$, $D'$ must contain a vertex from $v_r, x$ and $v_t$ because $x$ is a degree 2 vertex.
If $x\in D'$, we could replace $x$ by any one of $x_r$ and $x_t$ so that the resulting set remains a dominating set.
So assume that $D'\cap \{v_r, v_t\}\ne \emptyset$.
Let $D=D'\cup \{v_r,v_t\}$.
Then $|D|\le |D'|+1$ and $D$ is a dominating set of $G$.
It then follows that $\gamma(G)\le \gamma(G')+1\le  \lceil (n'+k')/4\rceil+1= \lceil (n-4+k-1)/4\rceil+1 \le  \lceil (n+k)/4\rceil$,
contradicting that the assumption that $G$ is a counterexample.
This completes the proof of Theorem~\ref{thm:main1}. 
\qed

\section{Hamiltonian plane triangulations}

Let $G$ be a Hamiltonian plane triangulation and let $H$ be a Hamilton cycle in $G$.
We can think of $H$ as bounding a triangulated inner subgraph $G_{int}$ and a triangulated outer subgraph
$G_{ext}$ such that $G_{int}\cap G_{ext}= H$.
Suppose $v\in V(G)$.
We denote by $\deg_{int} (v)$ (respectively, $\deg_{ext} (v)$) the degree of vertex $v$ in $G_{int}$ (resp. in $G_{ext}$).
A {\em 2-vertex} $v$ of $G$ is a vertex satisfying either $\deg_{int}(v)=2$ or $\deg_{ext}(v)=2$.
A {\em 2-chord} is a chord of the Hamilton cycle $H$ joining two vertices which lie at distance two on the cycle $H$.
In other words, a 2-chord is a chord joining the two neighbors of a 2-vertex.  A triangle containing exactly two edges of the Hamilton cycle 
is called a {\em type-2 triangle} in \cite{BOV}, which contains a 2-vertex and a 2-chord joining the two neighbors 
of the 2-vertex in the Hamilton cycle.

 
A Hamiltonian plane triangulation $G$ may have many Hamilton cycles.
If a Hamilton cycle $H$ is fixed, all 2-vertices and 2-chords of $H$ are fixed. 
The following result shows that a Hamiltonian plane triangulation without a dominating vertex must contain a Hamilton cycle $H$
which doesn't contain three consecutive 2-vertices (or three consecutive type-2 triangles $T_1, T_2$ and $T_3$ such that 
$T_1\cap T_{2}$ and $T_2\cap T_3$ are two incident edges of $H$). 

\begin{lemma}[Brinkmann, Ozeki and  Van Cleemput, \cite{BOV}]\label{lem2.1}  
Every Hamiltonian plane triangulation with $\gamma (G)\ge 2$ contains a Hamilton cycle containing no three consecutive 2-vertices.
\end{lemma}

So let $G$ be a plane triangulation which has a Hamilton cycle $H$ which in turn does not contain three consecutive 2-vertices.
A spanning subgraph of $G$ consisting of $H$ and all the 2-chords of $H$ was called an {\em $(H,A,B,O)$-graph} in \cite{PYZ}.
Lemma~\ref{lem2.1} shows that every Hamiltonian plane triangulation $G$ without a dominating vertex has an $(H,A,B,O)$-graph
as a spanning subgraph. 



\begin{thm}[\cite{PYZ}]\label{thm3.2}
Every $(H,A,B,O)$-graph on $n$ vertices with 
at least $(n+1)/2$ 2-chords has domination number at most
$\lceil 2n/7\rceil$.
\end{thm}

We are now prepared to prove Theorem~\ref{thm:main2}.  \medskip


\noindent{\bf Proof of Theorem~\ref{thm:main2}.}
Let $G$ be a Hamiltonian plane triangulation.
If $\gamma (G)=1$, the result is trivial.
So suppose that $\gamma (G)\ge 2$. 
By Lemma~\ref{lem2.1}, $G$ contains a Hamilton cycle $H$ which has no three consecutive 2-vertices.
Let $G_{int}$ be the subgraph of $G$ containing $H$ and everything inside of $H$
and let $G_{ext}$ be the subgraph of $G$ containing $H$ and everything outside of $H$. 

Let $K$ be the spanning subgraph consisting of $H$ together with all 2-chords of $H$. 
If $K$ has at least $(n+1)/2$ 2-chords, then $\gamma(K)\le   \lceil 2n/7\rceil$.
Since $K$ is spanning,
$\gamma(G)\le \gamma(K)\le  \lceil 2n/7\rceil\le 5n/16$ for all $n\ge 23$. 

So, in the following, assume that $K$ has at most $n/2$ 2-chords.
It follows that either $G_{int}$ or $G_{ext}$ has at most $n/4$ 2-chords.
Without loss of generality, assume that $G_{int}$ has
no more than $n/4$ 2-chords.
Then $G_{int}$ is a maximal outerplane graph with $t\le n/4$ vertices of
degree 2.  By Theorem~\ref{thm:outplanar},
$\gamma (G)\le (n+t)/4\le (n+ n/4)/4= 5n/16$.  
This completes the proof of the theorem. 
\qed
\medskip

\noindent{\bf Remark.} For Hamiltonian plane triangulations with at least seven, but less than twenty-three vertices,
the bound in Theorem~\ref{thm:main2} may be valid.
It may be possible to prove this or use a computer to check it.
However, the process would most likely be quite tedious and we do not see the necessity of doing so
since we do not expect the bound $5n/16$ to be tight for any infinite family of plane triangulations.


\end{document}